\documentclass[11pt]{amsart}
\usepackage{amssymb}
\newtheorem{theorem}{Theorem}[section]

\newtheorem{proposition}[theorem]{Proposition}

\newtheorem{lemma}[theorem]{Lemma}
\theoremstyle{remark}
\newtheorem{remark}[theorem]{Remark}
\newtheorem{definition}[theorem]{Definition}

\newtheorem{example}[theorem]{Example}
\newtheorem{examples}[theorem]{Examples}
\newcommand\be{\begin{equation}\label}
\newcommand\ee{\end{equation}}

\newcommand\G{\mathcal{G}}

\renewcommand{\O}{\mathcal{O}}
\newcommand{\Co}{\mathcal{C}}

\newcommand{\U}{\on{U}}

\newcommand{\N}{\mathbb{N}}
\newcommand{\R}{\mathbb{R}}
\newcommand{\C}{\mathbb{C}}

\newcommand{\Z}{\mathbb{Z}}

\newcommand\lie[1]{\mathfrak{#1}}
\renewcommand{\k}{\lie{k}}

\newcommand{\g}{\lie{g}}

\renewcommand{\t}{\lie{t}}
\newcommand{\Alc}{\lie{A}}

\newcommand{\on}{\operatorname}

\newcommand{\Ad}{ \on{Ad} }

\newcommand{\Hom}{ \on{Hom}} 
\renewcommand{\ker}{ \on{ker}} 
 
\newcommand{\SU}{ \on{SU}}

\newcommand\dirac{/\kern-1.2ex\partial} 
\newcommand\qu{/\kern-.7ex/} 

\newcommand{\lra}{\longrightarrow}
\newcommand{\hra}{\hookrightarrow}

\renewcommand{\d}{{\mbox{d}}}
\newcommand{\ol}{\overline}

\newcommand\sig{\sigma}

\newcommand\Om{\Omega}
\newcommand\om{\omega}

\newcommand{\f}{\frac}

\renewcommand{\l}{\langle}
\renewcommand{\r}{\rangle}
\newcommand{\hh}{{\textstyle \f{1}{2}}}
\newcommand{\ti}{\tilde}

\newcommand\pt{\on{pt}}

\newcommand\beqn{\begin{equation}}      
\newcommand\eeqn{\end{equation}}      
\newcommand{\ca}{\mathcal}
\newcommand{\wh}{\widehat}
\newcommand{\wt}{\widetilde}

\newcommand{\beq}{\begin{eqnarray*}}
\newcommand{\eeq}{\end{eqnarray*}}

\begin{document}

\title{The basic gerbe over a compact simple Lie group}
\author{Eckhard Meinrenken}
\address{University of Toronto, Department of Mathematics,
100 St George Street, Toronto, Ontario M5S3G3, Canada }
\email{mein@math.toronto.edu}
\begin{abstract}
Let $G$ be a compact, simply connected simple Lie group. We give a
construction of an equivariant gerbe with connection on $G$, with
equivariant 3-curvature representing a generator of
$H^3_G(G,\Z)$. Among the technical tools developed in this context 
is a gluing construction for equivariant bundle gerbes.
\end{abstract}
\maketitle

\section{Introduction}
Let $G$ be a compact, simply connected simple Lie group, acting on
itself by conjugation. It is well-known that the cohomology of $G$,
and also its equivariant cohomology, is trivial in degree less than
three and that $H^3(G,\Z)$ and $H^3_G(G,\Z)$ are canonically
isomorphic to $\Z$. The generator of $H^3(G,\Z)$ is represented by a
unique bi-invariant differential form $\eta\in \Om^3(G)$, admitting an
equivariantly closed extension $\eta_G\in\Om^3_G(G)$ in the complex of
equivariant differential forms.  Our goal in this paper is to give an
explicit, finite-dimensional description of an equivariant gerbe over
$G$, with equivariant 3-curvature $\eta_G$.

A number of constructions of gerbes over compact Lie groups may be
found in the literature, using different models of gerbes and valid in
various degrees of generality. The differential geometry of gerbes was
initiated by Brylinski's book \cite{br:lo}, building on earlier work
of Giraud. In this framework gerbes are viewed as sheafs of groupoids
satisfying certain axioms.  Brylinski gives a general construction of
a gerbe with connection, for any integral closed 3-form on any
2-connected manifold $M$.  The argument uses the path fibration
$P_0M\to M$, and is similar to the well-known construction of a line
bundle with connection out of a given integral closed 2-form on a
simply connected manifold.  In a later paper \cite{br:ge1}, Brylinski
gives a finite-dimensional description of the sheaf of groupoids
defining the basic gerbe for any compact Lie group $G$. A less
abstract picture, developed by Chatterjee-Hitchin
\cite{ch:co,hi:le,hi:wh}, describes gerbes in terms of {\em transition
line bundles} similar to the presentation of line bundles in terms of
transition functions.  A detailed construction of transition line
bundles for the basic gerbe over $G=\SU(N)$, (as well as for the much
more complicated case of finite quotients of $G=\SU(N)$) was obtained by
Gaw{\c{e}}dzki-Reis \cite{ga:wz}.

In this paper, we will extend the Gaw{\c{e}}dzki-Reis approach from
$\SU(N)$ to other simply connected simple Lie groups $G$. A
fundamental difficulty in the more general case is that, in contrast
to the case $G=\SU(N)$, the pull-back of a generator of $H^3_G(G,\Z)$
to a conjugacy class $\Co\subset G$ may not vanish. In this case it is
impossible to describe the basic gerbe in terms of a $G$-invariant
cover and $G$-equivariant transition line bundles. Compare with the
case of $G$-equivariant line bundles over $G$-manifolds $M$: Such a
line bundle may be described in terms of a $G$-invariant cover and
$G$-invariant transition functions only if its pull-back to any
$G$-orbit is equivariantly trivial.
 
One way of getting around this problem is to extend the
Chatterjee-Hitchin theory to the equivariant case, as in
\cite[Appendix A]{br:ge1}. A lift of the group action to a given gerbe
is obtained by specifying the isomorphisms between the gerbe and its
pull-back under the action of group elements $g\in G$. 
Unfortunately, the conditions for such isomorphisms to define a group 
action become rather complicated. A second possibility, 
adopted in this paper, is to use Murray's theory of 
{\em bundle gerbes} \cite{mu:bu}. 

To explain our approach in more detail, let us first discuss the
simplest case of $G=\SU(d+1)$, where it is equivalent to 
the construction in Gaw{\c{e}}dzki-Reis.  
The eigenvalues of any matrix $A\in \SU(d+1)$
can be uniquely written in the form
$$ \exp(2\pi i \lambda_1(A)),\ldots,\exp(2\pi i \lambda_{d+1}(A))$$ 
where $\lambda_1(A),\ldots,\lambda_{d+1}(A)\in\R$ satisfy
$\sum_{i=1}^{d+1}\lambda_i(A)=0$ and
$$ \lambda_1(A)\ge \lambda_2(A)\ge \cdots \ge \lambda_{d+1}(A) \ge 
\lambda_1(A)-1.$$
Define an open cover $V_1,\ldots,V_d,V_{d+1}$ of $G$, where $V_j$
consists of those matrices $A$ for which the $j$th inequality becomes
strict. Over the set $G_{\on{reg}}$ of regular elements, where all
inequalities are strict, we have $d+1$ line bundles
$L_1,\ldots,L_d,L_{d+1}$ defined by the eigenlines for the
eigenvalues $\exp(2\pi i \lambda_j(A))$.  For $i<j$, the tensor
product $L_{i+1}\otimes \cdots\otimes L_j\to G_{\on{reg}}$ extends to
a line bundle $L_{ij}\to V_i\cap V_j$. (One may view $L_{ij}$ as the
top exterior power of the sum of eigenspaces for the eigenvalues in
the given range.) For $i<j<k$ we have a canonical isomorphism
$L_{ij}\otimes L_{jk}\cong L_{ik}$ over the triple intersection
$V_i\cap V_j\cap V_k$.  The $L_{ij}$, together with these
isomorphisms, define a gerbe over $\SU(d+1)$, representing the
generator of $H^3(\SU(d+1),\Z)$.

More generally, consider any compact, simply connected, simple Lie
group $G$ of rank $d$.  Up to conjugacy, $G$ contains exactly $d+1$
elements with semi-simple centralizer. (For $G=\SU(d+1)$, these are
the central elements.) Let $\Co_1,\ldots,\Co_{d+1}\subset G$ be their
conjugacy classes.  We will define an invariant open cover
$V_1,\ldots,V_{d+1}$ of $G$, with the property that each member of this
cover admits an equivariant retraction onto the conjugacy class
$\Co_j\subset V_j$.  It turns out that every semi-simple
centralizer has a distinguished central extension by $\U(1)$. 
This central extension defines an equivariant bundle gerbe on 
$\Co_j$, hence (by pull-back) an equivariant bundle gerbe over $V_j$. 
We will find that these gerbes over $V_j$ glue together to produce 
a gerbe over $G$, using a gluing rule developed in this paper.

The organization of the paper is as follows. In Section
\ref{sec:gerbe} we review the theory of gerbes and pseudo-line bundles
with connections, and discuss 'strong equivariance' under a group
action.  Section \ref{sec:glue} describes gluing rules for bundle
gerbes. Section \ref{sec:prin} summarizes some facts about gerbes
coming from central extensions.  In Section \ref{sec:basic} we give
the construction of the basic gerbe over $G$ outlined above, and in
Section \ref{sec:conj} we study the ``pre-quantization of conjugacy
classes''. 

\vskip.1in
\noindent{\bf Acknowledgment:} I would like to thank Ping Xu for
fruitful discussions at the Poisson 2002 meeting in Lisbon, and for a
preliminary version of his preprint \cite{beh:eq} with Behrend and
Zhang, giving yet another construction of the basic gerbe over $G$.
Their (infinite-dimensional) approach is based on the notion of Morita
equivalence of (quasi-)symplectic groupoids. I thank the referees 
for detailed comments and suggestions. 
 
\section{Gerbes with connections}\label{sec:gerbe}
In this section we review gerbes on manifolds, along the lines of
Chatterjee-Hitchin and Murray. 

\subsection{Chatterjee-Hitchin gerbes}
Let $M$ be a manifold. Any Hermitian line bundle over $M$ 
can be described by an open cover $U_a$, and transition functions 
$\chi_{ab}:\,U_a\cap U_b\to \U(1)$ satisfying a cocycle condition 
$(\delta \chi)_{abc}=\chi_{bc}\chi_{ac}^{-1}\chi_{ab}=1$ on triple intersections. 
The cohomology class in $H^1(M,\underline{\U(1)})=H^2(M,\Z)$ 
defined by this cocycle is the Chern class of the line bundle. 
Chatterjee-Hitchin \cite{ch:co,hi:le,hi:ta} suggested to realize classes in $H^3(M,\Z)$ 
in a similar fashion, replacing $\U(1)$-valued functions with 
Hermitian line bundles. They define a gerbe to be a collection of 
Hermitian transition line bundles $L_{ab} \to U_a\cap U_b$ and and a 
trivialization, i.e. unit length section, $t_{abc}$ 
of the line bundle $(\delta L)_{abc}=L_{bc}L_{ac}^{-1}L_{ab}$ 
over triple intersections. These trivializations have to 
satisfy a compatibility relation over quadruple intersections, 
$$(\delta t)_{abcd}\equiv t_{bcd}t_{acd}^{-1}t_{abd}t_{abc}^{-1}=1,$$ 
which makes sense since $(\delta
t)_{abcd}$ is a section of the {\em canonically} trivial bundle. (Each
factor $L_{ab}$ cancels with a factor $L_{ab}^{-1}$.) After passing to
a refinement of the cover, such that all $L_{ab}$ become
trivializable, and picking trivializations, $t_{abc}$ is simply a
\v{C}ech cocycle of degree 2, hence defines a class in
$H^2(M,\underline{\U(1)})=H^3(M,\Z)$. The class is independent of 
the choices made in this construction, and is called the 
{\em Dixmier-Douady class} of the gerbe. 

Note that in practice, it is often not desirable to pass to a
refinement. For example, if $M$ is a connected, oriented 3-manifold,
the generator of $H^3(M,\Z)=\Z$ can be described in terms of the
cover $U_1,U_2$ where $U_1$ is an open ball around a given point $p\in
M$, and $U_2=M\backslash\{p\}$, using the degree one line bundle over
$U_1\cap U_2\cong S^2\times (0,1)$. 

\subsection{Bundle gerbes}\label{subsec:murray}
Bundle gerbes were invented by Murray \cite{mu:bu}, generalizing the following 
construction of line bundles. Let $\pi:\,X\to M$ be a fiber bundle, or 
more generally a surjective 
submersion. (Different components of $X$ may have different dimensions. ) 
For each $k\ge 0$ let $X^{[k]}$ denote the $k$-fold 
fiber product of $X$ with itself. There are $k+1$ projections 
$\partial^i:\,X^{[k+1]}\to X^{[k]}$, omitting the $i$th factor in the 
fiber product. Suppose we are given a smooth function 
$\chi:\,X^{[2]}\to \U(1)$, satisfying a cocycle condition 
$\delta\chi=1$ where 
$$\delta \chi:=\partial_0^*\chi\partial_1^*\chi^{-1}\partial_2^*\chi:\,X^{[3]}\to
\U(1).$$
Then $\chi$ determines a Hermitian line bundle $L\to M$, with fibers
at $m\in M$ the space of all linear maps 
$\phi:\,X_m=\pi^{-1}(m)\to \C$ such that 
$ \phi(x)=\chi(x,x')\phi(x')$. Given local sections $\sig_a:\,U_a\to X$ of 
$X$, the pull-backs of $\chi$ under the maps 
$(\sig_a,\sig_b):\,U_a\cap U_b\to X^{[2]}$ give transition functions 
$\chi_{ab}$ for the line bundle. 

Again, replacing $\U(1)$-valued functions by line bundles in this
construction, one obtains a model for gerbes: A bundle gerbe is given
by a line bundle $L\to X^{[2]}$ and a trivializing section $t$ of the
line bundle $\delta L= \partial_0^*L\otimes \partial_1^*L^{-1}
\otimes\partial_2^*L$ over $X^{[3]}$, satisfying a compatibility
condition $\delta t=1$ over $X^{[4]}$ (which makes sense since $\delta
t$ is a section of the canonically trivial bundle $\delta\delta L$).
Given local sections $\sig_a:\,U_a\to X$, one can pull these data back
under the maps $(\sig_a,\sig_b):\,U_a\cap U_b\to X^{[2]}$ and $(
\sig_a,\sig_b,\sig_c):\,U_a\cap U_b\cap U_c\to X^{[3]}$ to obtain a
Chatterjee-Hitchin gerbe. The Dixmier-Douady class of $(X,L,t)$ is by
definition the Dixmier-Douady class of this Chatterjee-Hitchin gerbe;
again this is independent of all choices.  The Dixmier-Douady class
behaves naturally under tensor product, pull-back and duals.

Notice that Chatterjee-Hitchin gerbes may be viewed as a special 
case of bundle gerbes, with $X$ the disjoint union of the 
sets $U_a$ in the given cover.

\begin{remark}
In his original paper \cite{mu:bu} Murray considered bundle gerbes
only for fiber bundles, but this was found too restrictive.  In
\cite{mu:bu1},\cite{st:ge} the weaker condition (called ``locally
split'') is used that every point $x\in M$ admits an open neighborhood
$U$ and a map $\sig:\,U\to X$ such that
$\pi\circ\sig=\on{id}$. However, this condition seems insufficient in
the smooth category, as the fiber product $X\times_M X$ need not be a
manifold unless $\pi$ is a submersion.
\end{remark}

\subsection{Simplicial gerbes}
Murray's construction fits naturally into a  wider context 
of {\em simplicial gerbes}. 
We refer to Mostow-Perchik's notes of lectures by R. Bott \cite{mo:no}
and to Dupont's paper \cite{du:si}
for a nice introduction to simplicial manifolds, and to Stevenson
\cite{st:ge} for their appearance in the gerbe context.

Recall that a 
{\em simplicial manifold} $M_\bullet$ is a sequence of manifolds
$(M_n)_{n=0}^\infty$, together with {\em face maps}
$\partial_i:\,M_n\to M_{n-1}$ for $i=0,\ldots,n$ satisfying relations
$\partial_i\circ\partial_j=\partial_{j-1}\circ\partial_i$ for $i<j$.
(The standard definition also involves {\em degeneracy maps} but these
need not concern us here.) 
The {\em (fat) geometric realization} of
$M_\bullet$ is the topological space $\|M\|=\coprod_{n=1}^\infty
\Delta^n \times M_n/ \sim $, where $\Delta^n$ is the $n$-simplex and
the relation is $(t,\partial_i(x))\sim (\partial^i(t),x)$, for
$\partial^i:\,\Delta^{n-1}\to \Delta^n$ the inclusion as the $i$th
face.  A (smooth) simplicial map between simplicial 
manifolds $M_\bullet,M_\bullet'$ is a collection of smooth maps 
$f_n:\,M_n\to M_n'$ intertwining the face maps; such a map 
induces a map between the geometric realizations.

\begin{examples}
\begin{enumerate}
\item 
If $S$ is any manifold, one can define a simplicial manifold
$E_\bullet S$ where $E_nS$ is the $n+1$-fold cartesian product 
of $S$, and $\partial_j$ omits the $j$th factor. It is known 
\cite{mo:no} that the geometric realization $||ES||$ of this 
simplicial manifold is contractible.  More generally, if 
$X\to M$ is a fiber bundle with fiber $S$, one can define 
a simplicial manifold $E_nX:=X^{[n+1]}$, with face maps 
as in Section \ref{subsec:murray}. 
The geometric realization $||EX||$ becomes a fiber bundle 
over $M$ with contractible fiber $||ES||$. 
\item \cite{mi:co,se:cl}
For any Lie group $G$ there is a simplicial manifold 
$B_nG=G^{n}$. The face maps $\partial_i$ for $0<i<n$ are  
$$\partial_i(g_1,\ldots,g_n)=(g_1,\ldots,g_ig_{i+1},\ldots,g_n),$$
while $\partial_0$ omits the  first component and 
$\partial_n$ the last component. The map $\pi_n:\,E_nG\to B_nG$ given by 
$\pi_n(k_0,\ldots,k_n)=(k_0 k_1^{-1},\ldots,k_{n-1}k_n^{-1})$ 
is simplicial, and the induced map on geometric realizations 
is a model for the classifying bundle $EG\to BG$. 
\item
\cite{se:cl,mo:no}
If $\ca{U}=\{U_a, a\in A\}$ is an open cover of $M$, one defines a 
simplicial manifold 
$$ \ca{U}_nM:=\coprod_{(a_0,\ldots,a_n)\in A_n} U_{a_0\ldots a_n}$$
where $A_n$ is the set of all sequences $(a_0,\ldots,a_n)$ such that 
$U_{a_0\ldots a_n}:=U_{a_0}\cap \cdots \cap U_{a_n}$ is non-empty. 
The face maps are induced by the inclusions, 
$$ \partial_i:\,U_{a_0\ldots a_n}\hra  U_{a_0\ldots \wh{a_i}\ldots a_n}.$$
One may view this as a special case of (a), with $X=\coprod_{a\in A}U_a$.  
It is known \cite[Theorem 7.3]{mo:no} that 
$||\ca{U}M||$ is homotopy equivalent to $M$. 
\item 
\cite{beh:eq} The definitions of $E_nG$ and $B_nG$ extend to Lie
groupoids $G$ over a base $S$. If $s,t:\,G\to S$ are the source and
target maps, one defines $E_nG$ as the $n+1$-fold fiber product of $G$
with respect to the target map $t$. The space $B_nG$ for $n\ge 1$ is
the set of all $(g_1,\ldots,g_n)\in G^n$ with $s(g_j)=t(g_{j-1})$,
while $B_0G=S$. The definition of the face maps $\partial_j:\,B_nG\to
B_{n-1}G$ is as before for $n>1$, while for $n=1$, $\partial_0=t$ and
$\partial_1=s$. We have a simplicial map $E_nG\to B_nG$ defined just
as in the group case. 
\end{enumerate}
\end{examples}

The bi-graded space of differential forms $\Om^\bullet(M_\bullet)$
carries two commuting differentials $\d,\delta$, where $\d$ is the de
Rham differential and $\delta:\, \Om^k(M_n)\to \Om^{k}(M_{n+1})$ is an
alternating sum, $\delta\alpha=\sum_{i=0}^{n+1} (-1)^i \partial_i^*
\alpha$. It is known \cite[Theorem 4.2, Theorem 4.5]{mo:no} that the
total cohomology of this double complex is the (singular) cohomology
of the geometric realization, with coefficients in $\R$.

We will use the $\delta$ notation in many similar situations: For
instance, given a Hermitian line bundle $L\to M_n$, we define a
Hermitian line bundle $\delta L\to M_{n+1}$ as a tensor product,
$$ \delta L=\partial_0^*L\otimes \partial_1^*L^{-1}\otimes\cdots 
\otimes \partial_{n+1}^*L^\pm.
$$
The line bundle $\delta(\delta L)\to M_{n+1}$ 
is canonically trivial, due to the relations between face maps. 
If $\sig$ is a unitary section (i.e. a trivialization) 
of $L$, one uses a similar formula to define a unitary section 
$\delta\sig$ of $\delta L$. Then $\delta(\delta\sig)=1$ 
(the identity section of the trivial line bundle $\delta(\delta L)$).
For any unitary connection $\nabla$ of $L$, one defines a unitary connection 
$\delta\nabla$ of $\delta L$ in the obvious way. 
{\bf For the rest of this paper, we take all line bundles $L$ to be {\em
Hermitian} line bundles, and all connections $\nabla$ on $L$ to be
{\em unitary} connections}.

Let $M_\bullet$ be a simplicial manifold. One might define a
simplicial line bundle as a collection of line bundles $L_n\to M_n$
such that the face maps $\partial_i:\,M_n\to M_{n-1}$ lift to line
bundle homomorphisms $\hat{\partial}_i:\,L_n\to L_{n-1}$, satisfying
the face map relations.  Thus $L_\bullet$ is itself a simplicial
manifold, and its geometric realization $\|L\|$ is a line bundle over
$\|M\|$. Equivalently, the lifts $\hat{\partial}_i$ may be viewed as
isomorphisms, $\partial_i^*L_{n-1}\to L_n$. In particular, we may
identify $L_n$ with the pull-back of $L:=L_0$ under the $n$th-fold
iterate $\partial_0\circ \cdots \circ \partial_0$.

The isomorphisms $\partial_1^*L\cong \partial_0^*L=L_1$ determine a
unitary section $t$ of $\delta L\to M_1$, and the compatibility of
isomorphisms
$$(\partial_0\partial_2)^*L\cong
(\partial_0\partial_1)^*L \cong (\partial_0\partial_0)^*L=L_2$$ 
amount to the condition $\delta t=1$. (Compatibility of the
isomorphisms for $L_n$ with $n\ge 3$ is then automatic.) That is, {\em
a simplicial line bundle over $M_\bullet$ is given by a line bundle
$L\to M_0$, together with a unitary section $t$ of $\delta L\to M_1$,
such that $\delta t=1$ over $M_2$.} A unitary section $s$ of $L$ with
$\delta s=t$ induces a unitary section of $\|L\|\to \|M\|$.

Taking $L$ to be trivial, we see in particular that any $\U(1)$-valued
function $t$ on $M_1$, with $\delta t=1$, defines a line bundle over the
geometric realization. A trivialization of that line bundle is given 
by a $\U(1)$-valued function on $M_0$ satisfying $\delta s=t$. 
Replacing $\U(1)$-valued functions with line bundles, 
this motivates the following definition.

\begin{definition}
A {\em simplicial gerbe} over $M_\bullet$ is a pair $(L,t)$,
consisting of a line bundle $L\to M_1$, together with a section $t$ of
$\delta L\to M_2$ satisfying $\delta t=1$. A pseudo-line bundle for
$(L,t)$ is a pair $(E,s)$, consisting of a line bundle $E\to M_0$ and
a section $s$ of $\delta E^{-1}\otimes L$ such that $\delta s=t$.
\end{definition}

\begin{remark}
\begin{enumerate}
\item 
We are using the notion of a simplicial gerbe only as a ``working
definition''. It is clear from the discussion above that a more general
notion would involve a gerbe over $M_0$.
\item
In \cite{br:ge1}, what we call simplicial gerbe is called a simplicial line 
bundle. The name pseudo-line bundle is adopted from 
\cite{br:ge1}, where it is used in a similar context.
\end{enumerate}
\end{remark}

A simplicial gerbe over $\ca{U}_\bullet M$ (for a cover $\ca{U}$ of
$M$) is a Chatterjee-Hitchin gerbe, while a simplicial gerbe over
$E_\bullet X=X^{[\bullet+1]}$ (for a surjective submersion $X\to M$)
is a bundle gerbe.  It is shown in \cite{mu:bu} that the
characteristic class of a bundle gerbe $(X,L,t)$ vanishes if and only if it
admits a pseudo-line bundle.

\begin{example}[Central extensions] (See \cite[p. 615]{br:ge1}.)
Let $K$ be a Lie group. A simplicial line bundle over $B_\bullet K$ is
the same thing as a group homomorphism $K\to \U(1)$: The line bundle 
$L\to B_0K$ is trivial since $B_0K$ is just a point, hence the unitary 
section $t$ of $\delta L$ becomes a $\U(1)$-valued function. The condition 
$\delta t=1$ means that this function is a group homomorphism. 

Similarly, a simplicial gerbe $(\Gamma,\tau)$ over $B_\bullet K$ is 
the same thing as a central extension 
$$\U(1)\to \wh{K}\to K$$
Indeed, given the line bundle $\Gamma\to K$ let $\wh{K}$ be 
the unit circle bundle inside $\Gamma$. The fiber of $\delta\Gamma\to K^2$ 
at $(k_1,k_2)$ is a tensor product $\Gamma_{k_2}\Gamma_{k_1k_2}^{-1}
\Gamma_{k_1}$, hence the section 
$\tau$ of $\delta\Gamma\to K^2$ defines a unitary isomorphism 
$\Gamma_{k_1}\Gamma_{k_2}\cong \Gamma_{k_1k_2}$, or equivalently 
a product on $\wh{K}$ covering the group 
multiplication on $K$. Finally, the condition  $\delta\tau=1$ is equivalent
to associativity of this product. 

A pseudo-line bundle $(E,s)$ for the simplicial
gerbe $(\Gamma,\tau)$ is the same thing as a splitting of the central
extension: Obviously $E$ is trivial since $B_0K$ is just a point; the
section $s$ defines a trivialization $\wh{K}=K\times\U(1)$, and
$\delta s=t$ means that this is a group homomorphism.
\end{example}

\begin{definition}
A connection on a simplicial gerbe $(L,t)$ over $M_\bullet$ 
is a line bundle connection $\nabla^L$, together with a 2-form 
$B\in \Om^2(M_0)$, such that  $(\delta\nabla^L)\,t=0$ and 
$$ \delta B=\f{1}{2\pi i}\on{curv}(\nabla^L).$$
Given a pseudo-line bundle $\ca{L}=(E,s)$, we say that $\nabla^E$ 
is a pseudo-line bundle connection if it has the property
$((\delta\nabla^E)^{-1}\nabla^L) s=0.$
\end{definition}

Simplicial gerbes need not admit connections in general. A sufficient 
condition for the existence of a connection is that the 
$\delta$-cohomology of the double complex $\Om^k(M_n)$ vanishes in 
bidegrees $(1,2)$ and $(2,1)$. In particular, this holds true 
for bundle gerbes: Indeed it is shown in \cite{mu:bu} 
that for any surjective submersion $\pi:\,X\to M$ the  sequence 
\begin{equation}\label{eq:exseq}
0\lra \Om^k(M)
\stackrel{\pi^*}{\lra}\Om^k(X)
\stackrel{\delta}{\lra}\Om^k(X^{[2]})
\stackrel{\delta}{\lra} \Om^k(X^{[3]})
\stackrel{\delta}{\lra} \cdots 
\end{equation}
is exact, so the $\delta$-cohomology vanishes in {\em all} degrees.

Thus, every bundle gerbe $\ca{G}=(X,L,t)$ over a manifold $M$ (and in
particular every Chatterjee-Hitchin gerbe) admits a connection.  One
defines the {\em 3-curvature} $\eta\in \Om^3(M)$ of the bundle gerbe
connection by $ \pi^*\eta=\d B\in\ker\delta$. It can be shown that its
cohomology class is the image of the Dixmier-Douady class $[\G]$ under
the map $H^3(M,\Z)\to H^3(M,\R)$. Similarly, if $\G$ admits a
pseudo-line bundle $\ca{L}=(E,s)$, one can always choose a pseudo-line
bundle connection $\nabla^E$. The difference $\f{1}{2\pi
i}\on{curv}(\nabla^E)-B$ is $\delta$-closed and one defines the {\em
error 2-form} of this connection by
$$ \pi^*\om=\f{1}{2\pi i}\on{curv}(\nabla^E)-B$$
It is clear from the definition that $\d\om +\eta=0$.

\begin{remark}
There is a notion of holonomy around surfaces for gerbe connections
(cf. Hitchin \cite{hi:le} and Murray \cite{mu:bu}), and 
in fact gerbe connections can be
defined in terms of their holonomy (see Mackaay-Picken \cite{mac:ho}).
\end{remark}

\subsection{Equivariant bundle gerbes}\label{subsec:strong}
Suppose $G$ is a Lie group acting on $X$ and on $M$, and that $\pi:\,X\to
M$ is a $G$-equivariant surjective submersion. Then $G$  acts on
all fiber products $X^{[p]}$. We will say that a bundle gerbe $\G=(X,L,t)$ is
{\em $G$-equivariant}, if $L$ is a $G$-equivariant line bundle
and $t$ is a $G$-invariant section. An equivariant bundle gerbe 
defines a gerbe over the Borel construction $X_G=EG\times_G X
\to M_G=EG\times_G M$ 
\footnote{We have not discussed bundle gerbes over
infinite-dimensional spaces such as $M_G$.  Recall however
\cite{bo:di} that the classifying bundle $EG\to BG$ may be
approximated by finite-dimensional principal bundles, and that
equivariant cohomology groups of a given degree may be computed using
such finite dimensional approximations.}, hence has an {\em
equivariant} Dixmier-Douady class in
$H^3(M_G,\Z)=H^3_G(M,\Z)$. Similarly, we say that a pseudo-line bundle
$(E,s)$ for $(X,L,t)$ is  equivariant, provided $E$ carries a
$G$-action and $s$ is an invariant section.

\begin{remark}
As pointed out in Mathai-Stevenson \cite{ma:ch}, this notion of
equivariant bundle gerbe is sometimes 'really too strong': For
instance, if $X=\coprod U_a$, for an open cover $\ca{U}=\{U_a,a\in
A\}$, a $G$-action on $X$ would amount to the cover being
$G$-invariant.  Brylinski \cite{br:ge1} on the other hand gives a
definition of equivariant Chatterjee-Hitchin gerbes that does not
require invariance of the cover. 
\end{remark}

To define  equivariant connections and curvature, we will need
some notions from equivariant de Rham theory \cite{gu:su}.  Recall
that for a compact group $G$, the equivariant cohomology
$H^\bullet_G(M,\R)$ may be computed from Cartan's complex of
equivariant differential forms $\Om_G^\bullet(M)$, consisting of
$G$-equivariant polynomial maps $\alpha:\,\g\to \Om(M)$. The grading
is the sum of the differential form degree and twice the polynomial
degree, and the differential reads
$$(\d_G\alpha)(\xi)=\d\alpha(\xi)-\iota(\xi_M)\alpha(\xi),$$
where
$\xi_M=\f{d}{d t}|_{t=0}\exp(-t\xi)$ is the generating vector field
corresponding to $\xi\in\g$. Given a $G$-equivariant connection
$\nabla^L$ on an equivariant line bundle, one defines \cite[Chapter
7]{be:he} a $\d_G$-closed equivariant curvature
$\on{curv}_G(\nabla^L)\in \Om^2_G(M)$.

A  equivariant connection on a  $G$-equivariant bundle
gerbe $(X,L,t)$ over $M$ is a pair $(\nabla^L,B_G)$, where $\nabla^L$
is an invariant connection and $B_G\in \Om^2_G(X)$ an equivariant
2-form, such that $\delta\nabla^L t=0$ and $\delta B_G=\f{1}{2\pi
i}\on{curv}_G(\nabla^L)$.  Its equivariant 3-curvature
$\eta_G\in\Om^3_G(M)$ is defined by $\pi^*\eta_G=\d_G B_G$. Given an
{\em invariant} pseudo-line bundle connection $\nabla^E$ on a
 equivariant pseudo-line bundle $(E,s)$, one defines the
equivariant error 2-form $\om_G$ by
$$\pi^*\om_G=\f{1}{2\pi i}\on{curv}_G(\nabla^E)-B_G.$$
Clearly, $\d_G\om_G+\eta_G=0$.

\section{Gerbes from principal bundles}
\label{sec:prin}
The following well-known example \cite{br:ge}, \cite{mu:bu} of a gerbe
will be important for our construction of the basic gerbe over
$G$. Suppose $\U(1)\to\wh{K}\to K$ is a central extension, and
$(\Gamma,\tau)$ the corresponding simplicial gerbe over
$B_\bullet K$. Given a principal $K$-bundle $\pi:\,P\to B$, one
constructs a bundle gerbe $(P,L,t)$, sometimes called the lifting
bundle gerbe. Observe that 
$$E_n P=P\times_K E_n K ,$$
which we may view as a fiber bundle over $B$ but also as a 
fiber bundle $E_nK\times_K P$ over $B_nK$.  Let
\begin{equation}\label{eq:f}
f_\bullet:\,E_\bullet P\to B_\bullet K\ee
be the bundle projection. Then $L=f_1^*\Gamma,\ \ t=f_2^*\tau$ defines a
bundle gerbe $(P,L,t)$. A pseudo-line bundle for this
bundle gerbe is equivalent to a lift of the structure group to
$\wh{K}$: Indeed if $\wh{P}$ is a principal $\wh{K}$-bundle lifting
$P$, consider the associated bundle $E=\wh{P}\times_{\U(1)}\C$.  From
the action map $\wh{K}\times \wh{P}\to \wh{P}$ one obtains an
isomorphism $\Gamma_k \otimes E_p\cong E_{k.p}$, or equivalently a
section $s$ of $\delta E^{-1}\otimes L$. One checks that $\delta s=t$,
so that $(E,s)$ is a pseudo-line bundle.  Conversely, the bundle
$\wh{P}$ is recovered as the unit circle bundle in $E$, and $s$
defines an action of $\wh{K}$ lifting the action of $K$.  See Gomi
\cite{go:con} for a detailed construction of bundle gerbe connections
on $(P,L,t)$.

\begin{remark}
To obtain a Chatterjee-Hitchin gerbe from this bundle gerbe, we must
choose a cover $\ca{U}$ of $M$ such that $P$ is trivial over each
$U_a\in\ca{U}$. Any choice of trivialization gives a simplicial map
$\ca{U}_\bullet M\to E_\bullet P$, and we pull back the bundle gerbe
under this map. More directly, the local trivializations give rise to a
'classifying map' $\chi_\bullet:\,\ca{U}_\bullet M\to B_\bullet K$
(see \cite{mo:no}), and the Chatterjee-Hitchin gerbe is defined as the
pull-back of $(\Gamma,\tau)$ under this map.
\end{remark}

Suppose the group $K$ is compact and connected.  After pulling back to
the universal cover $\wt{K}$, every central extension
$\U(1)\to\wh{K}\to K$ becomes trivial. It follows that every central
extension of $K$ by $\U(1)$ is of the form
$$\wh{K}=\wt{K}\times_{\pi_1(K)}\U(1),$$
where $\pi_1(K)\subset \ti{K}$
acts on $\U(1)$ via some homomorphism $\varrho\in
\on{Hom}(\pi_1(K),\U(1))$.  The choice of $\varrho$ for a given
extension is equivalent to the choice of a flat $\wh{K}$-invariant
connection on the principal $\U(1)$-bundle $\wh{K}\to K$. The central
extension is isomorphic to the {\em trivial} extension if and only if
$\varrho$ extends to a homomorphism $\wt{\varrho}:\,\wt{K}\to \U(1)$,
and the choice of any such $\wt{\varrho}$ is equivalent to a choice of
trivialization.  Using the natural map from
$(\k^*)^K=\on{Hom}(\wt{K},\R)$ onto $\on{Hom}(\wt{K},\U(1))$ this
gives an exact sequence of Abelian groups
\begin{equation}\label{eq:sequence}
(\k^*)^K\to \on{Hom}(\pi_1(K),\U(1))\to \{\mbox{central
extensions of }K\mbox{ by }\U(1)\} \to 1 .\end{equation}
Suppose $K$ is semi-simple (so that $(\k^*)^K=0$), and $T$ is a maximal torus
in $K$. Let $\wt{T}\subset \wt{K}$ be the maximal torus given as 
the pre-image of $T$. Let $\Lambda_K,\ti{\Lambda}_K\subset\t$ be the 
integral lattices of $T,\wt{T}$. The lattice $\ti{\Lambda}_K$ is equal to the 
co-root lattice of $K$, and  $\pi_1(K)=\Lambda_K/\ti{\Lambda}_K$,
(cf. \cite[Theorem V.7.1]{br:rep}). Therefore, if $K$ is semi-simple, 
$$ \{\mbox{central
extensions of }K\mbox{ by }\U(1)\} =\on{Hom}(\pi_1(K),\U(1))=\ti{\Lambda}^*_K/\Lambda^*_K,$$ 
the quotient of the dual of the co-root lattice by the weight lattice.

\begin{proposition}\label{prop:central}
Suppose $K$ is a compact, connected Lie group and $\pi:\,P\to M$ a principal 
$K$-bundle. 
\begin{enumerate}
\item
Any $\varrho\in\Hom(\pi_1(K),\U(1))$ defines a bundle gerbe $(P,L,t)$
over $M$, together with a gerbe connection $(\nabla^L,B)$ where $B=0$. In
particular this gerbe is {\em flat}.
\item
If $\varrho$ is the image of $\mu\in (\k^*)^K$, there is a
distinguished pseudo-line bundle $\ca{L}=(E,s)$ for this gerbe, with
$E$ a trivial line bundle.  Any principal connection $\theta\in\Om^1(P,\k)$
defines a connection on $\ca{L}$, with error 2-form $\om\in \Om^2(M)$ 
given by $\pi^*\om=\l\mu,F^\theta\r\in\Om^2(M)$, where $F^\theta$ is the
curvature.
\end{enumerate}
\end{proposition}
\begin{proof}
Let $\U(1)\to\wh{K}\to K$ be the central extension defined by
$\varrho$, and $(\Gamma,\tau)$ the corresponding simplicial gerbe over
$B_\bullet K$. As remarked above, $\varrho$ defines a flat connection
on $ \wh{K}\to K$, hence also a flat connection $\nabla^\Gamma$ on the
line bundle $\Gamma\to B_1 K$. Then $(\nabla^\Gamma,0)$ is a
connection on the simplicial gerbe $(\Gamma,\tau)$. Pulling back under
the map $f_\bullet$ (cf. \eqref{eq:f}) we obtain a connection
$(\nabla^L,0)$ on the bundle gerbe $(P,L,t)$.

If $\varrho$ is in the image of $\mu\in (\k^*)^K$, the corresponding
trivialization of $\wh{K}$ defines a unitary section $\sig$ of
$\Gamma$, with $\delta\sig=\tau$ and $\f{1}{2\pi
i}\nabla^\Gamma\sig=\l\mu,\theta^L\r\sig$, where $\theta^L$ is the
left-invariant Maurer-Cartan form on $K$.  Thus $\ca{L}=(E,s)$, with
$E$ the trivial line bundle and $s=f_1^*\sig$, is a pseudo-line
bundle for $\G$.
Given a principal connection $\theta$, let $\nabla^E$ be the
connection on the trivial bundle $E$, having connection 1-form
$\l\mu,\theta\r\in\Om^1(P)$.  
Since $\f{1}{2\pi i}\nabla^L s=f_1^*
\l\mu,\theta^L\r\,s$, it follows that 
\begin{equation}\label{eq:pse}
\f{1}{2\pi i} ((\delta\nabla^E)^{-1}\nabla^L) s=\l\mu,f_1^*\theta^L-\delta\theta\r.
\end{equation}
One finds $
\partial_1^*\theta=\Ad_{f_1^{-1}}(\partial_0^*\theta-f_1^*\theta^L)$.
Since $\mu$ is $K$-invariant, this shows that the right hand side of
\eqref{eq:pse} vanishes.  Thus $\nabla^E$ is a pseudo-line bundle
connection. The error 2-form $\om$ is given by
$$\pi^*\om=\d\l\mu,\theta\r=\l\mu,\d\theta\r=\l\mu,F^\theta\r.$$
\end{proof}

All of these constructions can be made equivariant in a rather 
obvious way: Thus if $G$ is another Lie group and $P$ is a
$G$-invariant principal $K$-bundle, any
$\varrho\in\Hom(\pi_1(K),\U(1))$ defines a  $G$-equivariant bundle
gerbe $(P,L,t)$ (with flat connection) over $M$.  If $\varrho$ is in
the image of $\mu\in (\k^*)^K$, there is a  $G$-equivariant
pseudo-line bundle for this gerbe. Furthermore any choice of
$G$-equivariant principal connection on $P$ defines a $G$-equivariant
pseudo-line bundle connection, with equivariant error 2-form
$\pi^*\om_G=\l\mu,F^\theta_G\r$ where $F^\theta_G\in \Om^2_G(P,\k)$ is
the equivariant curvature.

\section{Gluing data}\label{sec:glue}
In this Section we describe a procedure for gluing a collection of 
bundle gerbes $(X_i,L_i,t_i)$ on open subsets $V_i\subset M$, 
with pseudo-line bundles of their quotients on overlaps. 
\footnote{See Stevenson \cite{st:ge} for similar gluing constructions.} 
We begin with the somewhat simpler case that the surjective submersions 
$X_i\to V_i$ are obtained by restricting a surjective submersion 
$X\to M$, and later reduce the general case to this special case. 

Thus, let $\pi:\,X\to M$ be a surjective submersion and 
let $V_i,\ i=0,\ldots,d$ an open cover of $M$. Let $X_i=X|_{V_i}$, 
and more generally $X_I=X|_{V_I}$ where $V_I$ 
is the intersection of all $V_i$ with $i\in I$.

Suppose we are given bundle gerbes $(X_i,L_i,t_i)$ over $V_i$ and
pseudo-line bundles $(E_{ij},s_{ij})$ for the quotients 
$(X_{ij},L_j L_i^{-1},t_j t_i^{-1})$ over $V_i\cap V_j$, where
$E_{ij}=E_{ji}^{-1}$ and $s_{ij}=s_{ji}^{-1}$.  Note that
$E_{ij}E_{jk}E_{ki}$ is a pseudo-line bundle for the trivial gerbe,
hence is a pull-back $\pi^*F_{ijk}$ of a line bundle $F_{ijk}\to M$,
and we will also require a unitary section $u_{ijk}$ of that line
bundle.  Under suitable conditions the data $(E_{ij},s_{ij})$ and
$u_{ijk}$ can be used to 'glue' the gerbes $(X_i,L_i,t_i)$.  The glued
gerbe will be defined over the disjoint union $\coprod_{i=1}^d X_i$. 
We have
\beq \Big(\coprod_{i=1}^d X_i\Big)^{[2]}&=&\coprod_{ij} X_i\times_M X_j\\
     \Big(\coprod_{i=1}^d X_i\Big)^{[3]}&=
&\coprod_{ijk} X_i\times_M X_j\times_M X_k\\
\cdots&&
\eeq
Hence, the glued gerbe will be of the form $(\coprod_{i}
X_i,\coprod_{ij}L_{ij},\coprod_{ijk} t_{ijk})$ where $L_{ij}$ 
are line bundles over $X_i\times_M X_j$ and $t_{ijk}$ 
unitary sections of a line bundle $(\delta L)_{ijk}$ over 
$\coprod_{ijk} X_i\times_M X_j\times_M X_k$. We will define 
$L_{ij}$ by tensoring $L_i\to X^{[2]}$ (restricted to 
$X_i\times_M X_j$ with the pull-back of $E_{ij}$ 
under the map $\partial_1:\,X_i\times_M X_j\to X_{ij}$.

\begin{proposition}\label{prop:gerbeglue}
Suppose the sections $u_{ijk}$ satisfy the cocycle condition
$u_{jkl}u_{ikl}^{-1}u_{ijl} u_{ijk}^{-1}=1$, and the sections $s_{ij}$
satisfy a cocycle condition $s_{ij}s_{jk}s_{ki}=1$. Then there is a
well-defined gerbe $(\coprod_{i}
X_i,\coprod_{ij}L_{ij},\coprod_{ijk} t_{ijk})$ over $M$, where
$L_{ij}\to X_i\times_M X_j$ is the line bundle
$$ L_{ij}=L_j\otimes \partial_1^* E_{ij}$$ 
and 
$t_{ijk}$ is a section of $(\delta L)_{ijk}\to X_i\times_M X_j\times_M X_k$ 
given by 
\begin{equation}\label{eq:t}
 t_{ijk}=t_k\otimes\partial_2^*s_{kj} \otimes \partial_2^*\partial_1^* 
\pi^* u_{ijk}\end{equation}
\end{proposition}
\begin{proof}
A short calculation gives, 
$$ (\delta L)_{ijk}=(\delta L_k)\otimes \partial_2^*
(L_jL_k^{-1}\delta E_{kj}^{-1})\otimes \partial_2^*\partial_1^*
\pi^* F_{ijk}
$$
showing that $t_{ijk}$ is a well-defined section of $(\delta
L)_{ijk}$. One finds furthermore
\beq (\delta t)_{ijkl}&=&
(\delta t_l)\otimes \partial_3^*\Big(t_l t_k^{-1} \delta s_{kl}^{-1}
\otimes \partial_2^*\big( s_{lj}s_{jk}s_{kl}
\otimes \partial_1^*\pi^*(u_{jkl}u_{ikl}^{-1}
u_{ijl} u_{ijk}^{-1})\big)\Big)\\
&=&\partial_3^*\partial_2^*\big( s_{lj}s_{jk}s_{kl}
\otimes \partial_1^*\pi^*(u_{jkl}u_{ikl}^{-1}
u_{ijl} u_{ijk}^{-1})\big)
\eeq
which equals $1$ under the given assumptions on $u$ and $s$.
\end{proof}

The gluing construction described in this Proposition is
particularly natural for Chatterjee-Hitchin gerbes: Suppose $\ca{U}$
is an open cover of $M$, and $X=\coprod_{U\in\ca{U}}U$.  For any
decomposition $\ca{U}=\coprod_{i=1}^d\ca{U}_i$ let $V_i=\cup_{U\in
\ca{U}_i}U$, and $X_i=\coprod_{U\in\ca{U}_i}U$. Note that in this
case,
$$ \coprod_i X_i = X.$$
Suppose $(L_i,t_i)$ are Chatterjee-Hitchin gerbes for the cover
$\ca{U}_i$ of $V_i$, and that we are given pseudo-line bundles
$(E_{ij},s_{ij})$ and a section $u_{ijk}$ as above. Note that the
$E_{ij}$ are a collection of line bundles over intersections $U_a\cap
U_b$ where $U_a\in \ca{U}_i$ and $U_b\in \ca{U}_j$.  The gluing
construction gives a Chatterjee-Hitchin gerbe $(L,t)$ for the cover
$\ca{U}$ of $M$, where the $E_{ij}$ enter the definition of transition
line bundles between open sets in distinct $\ca{U}_i,\ca{U}_j$.

\begin{remark}
Suppose $X=M$, and that all $L_i,t_i,s_{ij}$ are trivial.
Then the gerbe described in Proposition \ref{prop:gerbeglue} {\em is}
a Chatterjee-Hitchin gerbe for the cover $\{V_i\}$. The $E_{ij}$ now play
the role of transition line bundles, and $u_{ijk}$ play the role of $t$.
\end{remark}

Suppose now that, in addition to the assumptions of Proposition
\ref{prop:gerbeglue}, we have gerbe connections
$(\nabla^{L_i},B_i)$ and pseudo-line bundle connections
$\nabla^{E_{ij}}=(\nabla^{E_{ji}})^{-1}$. Let $\om_{ij}$ denote the
error 2-form for $\nabla^{E_{ij}}$.

\begin{proposition}\label{prop:connglue}
The connections 
$ \nabla^{L_{ij}}=\nabla^{L_j}\otimes\partial_1^*\nabla^{E_{ij}}$
on $L_{ij}$, together with the two forms $B_i\in\Om^2(X_i)$, 
define a gerbe connection if all error 2-forms $\om_{ij}$ vanish, and 
if 
$$\nabla^{E_{ij}}\nabla^{E_{jk}}\nabla^{E_{ki}}(\pi^* u_{ijk})=0.$$ 
\end{proposition}
\begin{proof}
Let $B$ be the 2-form on $\coprod X_i$ given by $B_i$ on $X_i$. 
We first verify that $ \f{1}{2\pi i}\on{curv}(\nabla^{L_{ij}})=(\delta
B)_{ij}$: 
\beq \f{1}{2\pi i}\on{curv}(\nabla^{L_{ij}})
&=&\f{1}{2\pi i}\on{curv}(\nabla^{L_j})
+\f{1}{2\pi i}\partial_1^*\on{curv}(\nabla^{E_{ij}})\\
&=&\delta B_j+\partial_1^*(B^j-B^i+\pi^*\om_{ij})\\
&=&\partial_0^* B_j-\partial_1^*B_i=(\delta B)_{ij}.
\eeq
Next, we check that $t_{ijk}$ is parallel for $(\delta\nabla^L)_{ijk}$:
\beq 
(\delta\nabla^L)_{ijk}
&=&\partial_0^*\nabla^{L_{jk}}\partial_1^*(\nabla^{L_{ik}})^{-1}
\partial_2^*\nabla^{L_{ij}}\\
&=&\delta\nabla^{L_k}\otimes 
\partial_2^*(\nabla^{L_k}(\nabla^{L_j})^{-1}\delta \nabla^{E_{jk}})
\otimes
\partial_2^*\partial_1^*(\nabla^{E_{ij}}\nabla^{E_{jk}}\nabla^{E_{ki}}).
\eeq
This annihilates \eqref{eq:t} as required.
\end{proof} 

We now describe a slightly more complicated gluing construction, 
in which the $X_i$ are not simply the restrictions of a surjective submersion 
$X\to M$. Instead, we assume that for each $I$ we are given a 
surjective submersion $\pi_I:\,X_I\to V_I$ are
surjective submersions, and for each $I\supset J$ a
fiber preserving smooth map $f_I^J:\,X_I\to X_J$, with the 
compatibility condition $f_J^K\circ
f_I^J=f_I^K$ for $I\supset J\supset K$. Our gluing data will consist
of the following:
\begin{enumerate}
\item[(i)] Over each $V_i$, bundle gerbes 
$(X_i,L_i,t_i)$ with connections $(\nabla^{L_i},B_i)$,
\item[(ii)] Over each $V_{ij}$, pseudo-line bundles
$E_{ij}=E_{ji}^{-1},s_{ij}=s_{ji}^{-1}$ with connections
$\nabla^{E_{ij}}=(\nabla^{E_{ji}})^{-1}$ for the bundle gerbe
$(X_{ij},L_{ij},t_{ij})$, given as the quotient of the pull-back of
$(X_j,L_j,t_j)$ by $f^j_{ij}$ and the pull-back of $(X_i,L_i,t_i)$ by
$f^i_{ij}$.
\item[(iii)] Over triple intersections, unitary sections $u_{ijk}$ of
the line bundle $F_{ijk}\to V_{ijk}$ defined by tensoring the
pull-backs of $E_{ij},E_{jk},E_{ki}$ by the maps
$f^{ij}_{ijk},f^{jk}_{ijk},f^{ki}_{ijk}$.
\end{enumerate}
We require that the $s_{ij}$ and $u_{ijk}$ satisfy a cocycle condition
similar to Proposition \ref{prop:gerbeglue}, that all error 
2-forms $\om_{ij}$ are zero, and that the connections 
$\nabla^{E_{ij}}$ satisfy a compatibility condition as in 
\ref{prop:connglue}.

These data may be used to define a bundle gerbe over $M$, by reducing 
to the setting of Propositions \ref{prop:gerbeglue}, \ref{prop:connglue}.
As a first step we construct a more convenient cover. 

\begin{lemma}\label{lem:covers}
There are open subsets $U_I$ of $M$, with $\ol{U_I}\subset V_I$,
and $\bigcup_I U_I=M$, such
that 
$$\ol{U_I}\cap \ol{U_J}=\emptyset\mbox{
unless }J\subset I\mbox{ or }I\subset J.$$  
The collection of open subsets
$$ V_i'=M\backslash \bigcup_{J\not\ni i}\ol{U_J}$$
is a shrinking of the open cover $V_i$, 
that is, $\bigcup V_i'=M$ and $\ol{V'_i}\subset V_i$.
\end{lemma}

The proof of this technical Lemma is deferred to Appendix \ref{sec:cover}.
Now set $X=\coprod_I X_I|_{U_I}$. 
By definition of $V_i'$, the restriction $X_i'=X|_{V_i'}$
is given by 
$$ X_i'=\coprod_{J\ni i} X_J|_{U_J\cap V_i'}.$$
More generally, letting $V_I'=\bigcap_{i\in I}V_i'$ and 
$X_I'=X|_{V_I'}$ we have 
$$ X_I'=\coprod_{J\supset I} X_J|_{U_J\cap V_I'}.$$
Let $X_I'\to X_I|_{V_I'}$ be the fiber preserving map, given 
on $X_J|_{U_J\cap V_I'}$ by the map $f_J^I:\,X_J\to X_I$.
Using these maps, we can pull-back our gluing
data: Let $(X_i',L_i',t_i')$ be the pull-back of the bundle gerbe
$(X_i,L_i,t_i)$ under the map $X_i'\to X_i$, equipped with the
pull-back connection. On overlaps $V_{ij}'$, we let
$(E_{ij}',s_{ij}')$ be the pseudo-line bundle with connections defined
by pulling back $(E_{ij},s_{ij})$. The gluing data obtained in this
way satisfy the conditions from Propositions \ref{prop:gerbeglue} and
\ref{prop:connglue}, and hence give rise to a bundle gerbe with
connection over $M$.

\begin{remark}
In our applications, the line bundles $E_{ij}$ are in fact trivial, so
one can simply take $u_{ijk}=1$ in terms of the trivialization. The
$s_{ij}$ are $\U(1)$-valued functions in this case, and the
compatibility condition reads $s_{ij}s_{jk}s_{ki}=1$ over $X_{ijk}$.
\end{remark}

The gluing constructions generalize  equivariant bundle gerbes
in a straightforward way.

\section{The basic gerbe over a compact simple Lie group}
\label{sec:basic}
In this section we explain our construction of the basic gerbe over a
compact, simple, simply connected Lie group. 

\subsection{Notation}
Let $G$ be a compact, simple, simply connected Lie group, with Lie
algebra $\g$. For any action of $G\times M\to M,\ (g,m)\mapsto g.m$ on
a manifold $M$, we will denote by $G_m$ the stabilizer group of a
point $m\in M$. If $M=G$ or $M=\g$, we will always consider the
adjoint action of $G$ unless specified otherwise.  
For instance, $G_g$ for denotes the centralizer of an element 
$g\in G$.

Choose a maximal torus $T$ of $G$, with Lie algebra $\t$. Let
$\Lambda=\ker(\exp|_\t)$ be the integral lattice and
$\Lambda^*\subset\t^*$ its dual, the (real) weight lattice.
Equivalently, $\Lambda$ is characterized as the lattice generated by
the coroots $\check{\alpha}$ for the (real) roots $\alpha$.  Recall
that the {\em basic inner product} $\cdot$ on $\g$ is the unique
invariant inner product such that
$\check{\alpha}\cdot\check{\alpha}=2$ for all long roots $\alpha$.
Throughout this paper, we will use the basic inner product to identify
$\g^*\cong\g$. Choose a collection of simple roots
$\alpha_1,\ldots,\alpha_d\in\Lambda^*$ and let
$\t_+=\{\xi|\alpha_j\cdot \xi\ge 0,\ j=1,\ldots,d\}$ be the
corresponding positive Weyl chamber. The fundamental alcove 
$\Alc$ is the subset cut out from $\t_+$ by the additional 
inequality $\alpha_0\cdot\xi\ge -1$ where $\alpha_0$ is the lowest 
root. 

The fundamental alcove parametrizes conjugacy classes in $G$, in the 
sense that each conjugacy class contains a unique point $\exp\xi$ 
with $\xi\in\Alc$. The quotient map will be denoted 
$q:\,G\to \Alc$.  
Let $\mu_0,\ldots,\mu_d$ be the vertices of $\Alc$, with $\mu_0=0$. 
For any $I\subseteq\{0,\ldots,d\}$, all group elements $\exp\xi$ 
with $\xi$ in the open face spanned by $\mu_j$ with $j\in I$
have the same centralizer, denoted $G_I$. In particular, 
$G_j$ will denote the centralizer of ${\exp\mu_j}$. 

For each $j$ let $\Alc_j\subset \Alc$ be the open star at $\mu_j$,
i.e. the union of all open faces containing $\mu_j$ in their
closure. Put differently, $\Alc_j$ is the complement of the closed
face opposite to the vertex $\mu_j$. We will work with the open cover
of $G$ given by the pre-images, $V_j=q^{-1}(\Alc_j)$.  More generally
let $\Alc_I=\cap_{j\in I}\Alc_j$, and $V_I:=q^{-1}(\Alc_I)$. The
flow-out $S_I=G_I.\exp(\Alc_I)$ of $\exp(\Alc_I)\subset T$ under 
the action of $G_I$ is an open subset of $G_I$, and 
is a slice for the conjugation action of $G$. That is,
$$ G\times_{G_I}S_I=V_I.$$
We let $\pi_I:\,V_I\to G/G_I$ denote the projection to the base. 

\subsection{The basic 3-form on $G$}\label{subsec:basic}
Let $\theta^L,\theta^R\in\Om^1(G,\g)$ be the left- and right-invariant 
Maurer-Cartan forms on $G$, respectively. The
3-form $\eta\in\Om^3(G)$ given by 
\footnote{For $\g$-valued forms $\beta_1,\beta_2$, the
bracket $[\beta_1,\beta_2]$ denotes the $\g$-valued form obtained by
applying the Lie bracket $[\cdot,\cdot]:\,\g\otimes\g\to\g$ to
the $\g\otimes\g$-valued form $\beta_1\wedge\beta_2$.}
$$ \eta=\f{1}{12} \theta^L\cdot [\theta^L,\theta^L]
=\f{1}{12} \theta^R\cdot [\theta^R,\theta^R]$$
is closed, and has a closed
equivariant extension $\eta_G\in \Om^3_G(G)$ given by
$$ \eta_G(\xi):=\eta-\hh (\theta^L+\theta^R)\cdot \xi.$$
Their cohomology classes represent generators of $H^3(G,\Z)=\Z$ and
$H^3_G(G,\Z)=\Z$, respectively.  The pull-back of $\eta_G$ to any
conjugacy class $\iota_\Co:\,\Co\hra G$ is exact.  In fact, let 
$\om_\Co\in \Om^2(\Co)^G\subset\Om^2_G(\Co)$ be the invariant 
2-form given on generating vector fields 
$ \xi_\Co,\xi'_\Co$ for $\xi,\xi'\in\g$ 
by the formula 
$$ \om_\Co(\xi_\Co(g),\xi'_\Co(g))=\hh \xi\cdot (\Ad_g-\Ad_{g^{-1}})\xi'.$$
Then \cite{al:mom,gu:gr}
$$ \d_G\om_\Co +\iota_\Co^*\eta_G=0.$$
We will
now show that $\eta_G$ is exact over each of the open subsets
$V_j$. Let $\Co_j=q^{-1}(\mu_j)\subset V_j$ be the conjugacy classes
corresponding to the vertices.

\begin{lemma}
The linear retraction 
$$ [0,1]\times\Alc_j\to \Alc_j,\ (t,\mu_j+\zeta)\mapsto \mu_j+(1-t)\zeta$$
of $\Alc_j$ onto the vertex $\mu_j$ lifts uniquely to a smooth
$G$-equivariant retraction from $V_j$ onto $\Co_j$.
\end{lemma}
\begin{proof}
Recall that the slice $S_j$ is an open neighborhood of $\exp(\mu_j)$
in $G_j$. Any $G_j$-equivariant retraction from $S_j$ onto $\exp\mu_j$
extends uniquely to a $G$-equivariant retraction from
$V_j=G\times_{G_j}S_j$ onto $\Co_j$.  Note that
$S_j'=G_j.(\Alc_j-\mu_j)$ is a star-shaped open neighborhood of $0$ in
$\g_j$, and that $S_j'\to S_j,\,\zeta\mapsto \exp(\mu_j)\exp(\zeta)$
is a $G_j$-equivariant diffeomorphism. The linear retraction of $S_j'$
onto the origin gives the desired retraction of $S_j$.  Uniqueness is
clear, since the retraction has to preserve $\exp(\Alc_j)\subset V_j$,
by equivariance.
\end{proof}

Let
$${\bf h}_j:\,\Om^p(V_j)\to \Om^p([0,1]\times V_j)\to \Om^{p-1}(V_j)$$ 
be the deRham homotopy operator for this retraction, given (up to a
sign) by pull-back under the retraction, followed by integration over
the fibers of $[0,1]\times V_j\to V_j$. It has the property
\begin{equation}\label{eq:hom}
\d_G {\bf h}_j+ {\bf h}_j \d_G =\on{Id}-\pi_j^*\iota_j^*
\end{equation}
where $\iota_j:\,\Co_j\to V_j$ is the inclusion and 
$\pi_j:\,V_j=G\times _{G_j} S_j\to G/G_j=\Co_j$ the projection.
Let $(\varpi_j)_G={\bf h}_j\eta_G-\pi_j^*\om_{\Co_j}\in\Om^2_G(V_j)$, 
and write 
$ (\varpi_j)_G=\varpi_j-\Psi_j$
where $\varpi_j\in \Om^2(V_j)$ and $\Psi_j\in \Om^0(V_j,\g)$. 

\begin{proposition}\label{prop:geom}
The equivariant 2-form $ (\varpi_j)_G=\varpi_j-\Psi_j$ has the following 
properties. 
\begin{enumerate}
\item
$ \d_G(\varpi_j)_G=\eta_G$.
\item
The pull-back of $(\varpi_j)_G$ to a conjugacy class $\Co\subset V_j$
is given by 
$$\iota_\Co^*(\varpi_j)_G=\Psi_j^*(\om_\O)_G-\om_\Co,$$
where $(\om_\O)_G$ is the equivariant symplectic form on the
adjoint orbit $\O=\Psi_j(\Co)$,
\item 
The pull-back of $\Psi_j$ to the conjugacy class $\Co_j$ vanishes. 
In fact, $\Psi_j(\exp\xi)=\xi-\mu_j$ for all $\xi\in \Alc_j$. 
\item
Over each intersection $V_{ij}=V_i\cap V_j$, the difference $\Psi_i-\Psi_j$
takes values in the adjoint orbit $\O_{ij}$ through 
$\mu_j-\mu_i\in\g\cong\g^*$. Furthermore, 
$$ (\varpi_j)_G-(\varpi_i)_G=-p_{ij}^*(\om_{\O_{ij}})_G$$
where $p_{ij}:\,V_{ij}\to \O_{ij}$ is the map defined by $\Psi_i-\Psi_j$, 
and $(\om_{\O_{ij}})_G$ is the equivariant symplectic form on the orbit.
\end{enumerate}
\end{proposition}

\begin{proof}
(a) holds by construction. (b) follows from the observation that
$\iota_\Co^*(\varpi_j)_G+\om_\Co$ is an equivariantly closed 2-form on
$\Co_j$, with $\Psi_j$ as its moment map.  To prove (c) we note that
since the retraction is equivariant, we have $\ti{\bf h}_j\circ
(\exp|_{\Alc_j})^*=(\exp|_{\Alc_j})^*\circ {\bf h}_j$ where
$(\exp|_{\Alc_j})^*$ is pull-back to
$\Alc_j\subset\t$ and where $\ti{\bf h}_j$ is the homotopy operator
for the linear retraction of $\t$ onto $\{\mu_j\}$. Let
$\nu:\,{\Alc_j}\to\t$ be the coordinate function (inclusion). Then
$$\ti{\bf h}_j\circ (\exp|_{\Alc_j})^*\hh (\theta^L+\theta^R)=\ti{\bf h}_j\circ\d\nu= \nu-\mu_j,$$
proving that $(\exp|_{\Alc_j})^*\Psi_j=\nu-\mu_j$. This yields 
(c), by equivariance.
For $\nu\in\Alc_{ij}$ we have, using (c),  
$$ (\Psi_i-\Psi_j)(\exp\nu)=(\nu-\mu_i)-(\nu-\mu_j)=\mu_j-\mu_i.$$
By equivariance, it follows that $\Psi_i-\Psi_j$ takes values in the
adjoint orbit through $\mu_j-\mu_i$. The difference
$\varpi_i-\varpi_j$ vanishes on the maximal torus $T$,  
and is therefore determined by its
contractions with generating vector fields. Since $\Psi_i-\Psi_j$ is a
moment map for $\varpi_i-\varpi_j$, it follows that
$\varpi_i-\varpi_j$ equals the pull-back of the symplectic form on
$G.(\mu_j-\mu_i)$.
\end{proof}

\subsection{The special unitary group}
For the special unitary group $G=\SU(d+1)$, the construction of the
basic gerbe simplifies due to the fact that in this case all vertices
$\mu_j$ of the alcove are contained in the weight lattice. In fact the gerbe 
is presented as a Chatterjee-Hitchin gerbe for the cover 
$\ca{V}=\{V_i,\ i=0,\ldots,d\}$. 

For each weight $\mu\in\Lambda^*\subset \t\subset \g$, let $G_\mu$ be
its stabilizer for the adjoint action and let $\C_\mu$ the
1-dimensional $G_\mu$-representation with infinitesimal character
$\mu$. Let the line bundle $L_\mu=G\times_{G_\mu}\C_\mu$ equipped with
the unique left-invariant connection $\nabla$. Then $L_\mu$ is a
$G$-equivariant pre-quantum line bundle for the orbit $\O=G.\mu$. That
is,
$$\f{i}{2\pi} \on{curv}_G(\nabla)=(\om_{\O})_G:=\om_{\O}-\Phi_{\O}$$
where $\om_{\O}$ is the symplectic form and
$\Phi_{\O}:\,\O\hra \g^*$ is the moment map given as
inclusion.

In particular, in the case of $\SU(d+1)$ all orbits
$\O_{ij}=G.(\mu_j-\mu_i)$ carry $G$-equivariant pre-quantum line
bundles. Recall the fibrations $p_{ij}:\,V_{ij}\to \O_{ij}$ defined by
$\Psi_i-\Psi_j$, and let
$$ L_{ij}=p_{ij}^*\,(L_{\mu_j-\mu_i}),$$
equipped with the pull-back connection.  For any triple intersection
$V_{ijk} =G\times_{G_{ijk}}S_{ijk}$, the tensor product $(\delta
L)_{ijk}=L_{jk}L_{ik}^{-1}L_{ij}$ is the pull-back of the line bundle
over $G/G_{ijk}$, defined by the zero weight
$$ (\mu_k-\mu_j)-(\mu_k-\mu_i)+(\mu_j-\mu_i)=0
$$
of $G_{ijk}$. It is hence canonically trivial, with $(\delta\nabla)_{ijk}$
the trivial connection. The trivializing section $t_{ijk}=1$ satisfies 
$\delta t=1$ and $(\delta\nabla)t=0$. Take $(B_j)_G=(\varpi_j)_G$. 
Then 
$$ (B_j)_G-(B_i)_G=(\varpi_j)_G-(\varpi_i)_G=-p_{ij}^*(\om_{\O_{ij}})_G=
\f{1}{2\pi i}\on{curv}_G(\nabla^{L_{ij}}).$$
Thus $\ca{G}=(\ca{V},L,t)$ is a  equivariant gerbe with connection 
$(\nabla,B)$. Since 
$$\d_G (B_j)_G=\d_G(\varpi_j)_G=\eta_G|_{V_j},$$
this is the basic gerbe for $\SU(d+1)$. The transition line bundles $L_{ij}$ 
may be expressed in terms of eigenspace line bundles, leading to the
description of the basic gerbe from the introduction.  

\begin{remark}
This description of the basic gerbe over the special unitary group was
independently found by Gaw{\c{e}}dzki-Reis \cite{ga:wz}, who also
discuss the much more difficult case of quotients of $\SU(d+1)$ by
subgroups of the center.
\end{remark}

A similar construction works for the group $C_d=\on{Sp}(d)$, the only case
besides $A_d=\SU(d+1)$ for which the vertices of the alcove are in the
weight lattice.  The following table lists, for all simply connected
compact simple groups, the smallest integer $k_0>0$ such that
$k_0\Alc$ is a weight lattice polytope.\footnote{This information is
extracted from the tables in Bourbaki \cite{bo:li}. Letting 
$w_1,\ldots,w_d$ be the fundamental weights, one determines 
$k_0$ as the least common multiple of the numbers 
${\alpha_{max}\cdot w_j}$, using the basic inner product
defined by $\alpha_{max}\cdot\alpha_{max}=2$.}  The construction
for $\SU(d+1)$ generalizes to describe the $k_0$'th power of the basic
gerbe in all cases.
\begin{equation}\label{eq:table}
 \begin{array}{c||c|c|c|c|c|c|c|c|c} 
G& A_d & B_d & C_d & D_d & E_6 & E_7 & E_8 & F_4 & G_2 \\ 
\hline \hline\rule{0mm}{5mm}
k_0& 1&2&1&2&3&12&60&6&2 \\ 
\end{array}
\end{equation}

\subsection{The basic gerbe for general simple, simply connected $G$}
\label{subsec:central}
The extra difficulty for the groups with $k_0>1$ comes from the fact
that the pull-back maps $H^3_G(G,\Z)\to H^3_G(\Co_j,\Z)\cong 
H^3_G(V_j,\Z)$ may be a non-zero torsion class, in general. 
In this case the restriction of the basic gerbe to
$V_j$ will be non-trivial. Our strategy for the general case is to
first construct equivariant bundle gerbes over $V_j$, and then glue
the local data as explained in Section \ref{sec:glue}.

The centralizers $G_g$ of elements $g\in G$ are always
connected \cite[Corollary (3.15)]{du:li} but need not be
simply-connected.  The conjugacy classes $\Co_j=q^{-1}(\mu_j)$
corresponding to the vertices of the alcove are exactly the conjugacy
classes of elements for which the centralizer is semi-simple.
Since 
$$H^3_G(\Co_j,\Z)=H^3_G(G/G_j,\Z)=H^3_{G_j}(\pt,\Z),$$ 
we see that the torsion problem described above is related to a
possibly non-trivial central extension of the centralizers $G_j$ of
$\exp(\mu_j)$ by the circle $\U(1)$. 

\begin{proposition}
Any vertex $\mu_j$ of the alcove $\Alc$ is in the dual of the 
co-root lattice 
for the corresponding centralizer $G_j$. It hence defines a 
homomorphism $\varrho_j\in \Hom(\pi_1(G_j),\U(1))$, or equivalently 
a central extension of $G_j$ by $\U(1)$. 
\end{proposition}
\begin{proof}
Let $\wt{G_j}$ be the universal cover of $G_j$.  A system of simple
roots for $\wt{G_j}$ is given by the list of all $\alpha_i$
($i=0,\ldots,d$) with $j\not=i$.  The lattice $\Lambda_j$ is spanned
by the corresponding coroots $\check{\alpha}_i$. To show that $\mu_j$
is in the dual of the co-root lattice, we have to verify that
$\l\mu_j,\check{\alpha}_i\r\in \Z$ for $i\not=j$. For $i\not=0,j$ this
is obvious since $\mu_j(\check{\alpha}_i)=0$. For $i=0$, we have
$||\check{\alpha}_0||^2=2$, and therefore $\check{\alpha}_0=\alpha_0$
and $\mu_j(\check{\alpha}_0)=\alpha_0(\mu_j)=-1$.
\end{proof}
Recall that for $i\not=j$,
$G_{ij}$ is the centralizer of points $\exp\mu$ with  
$\mu=t\mu_j+(1-t)\mu_i$ for some $0<t<1$. 
Let $\varrho_{ij}\in\Hom(\pi_1(G_{ij}),\U(1))$ be the 
quotient of $\pi_1(G_{ij})\to \pi_1(G_j)\stackrel{\varrho_j}{\lra}\U(1)$ 
by the homomorphism $\pi_1(G_{ij})\to \pi_1(G_i)
\stackrel{\varrho_i}{\lra}\U(1)$. 
 
\begin{lemma}\label{lem:diff}
The difference $\mu_j-\mu_i\in \g_{ij}$ is fixed under $G_{ij}$, and
$\varrho_{ij}\in\Hom(\pi_1(G_{ij}),\U(1))$ is its image under the
exact sequence \eqref{eq:sequence} for $K=G_{ij}$.
\end{lemma}

\begin{proof}
Since $G_{ij}$ fixes the curve
$g(t)=\exp(t\mu_j+(1-t)\mu_i)=\exp(\mu_i)\exp(t(\mu_j-\mu_i))$, it
stabilizes the Lie algebra element $\mu_j-\mu_i$.  The second claim is
immediate from the definition.
\end{proof}

We are now in position to explain our construction 
of the basic gerbe in the general case. For all 
$I\subset\{0,\ldots,d\}$ let $X_I\to V_I$ be the $G$-equivariant 
principal $G_I$-bundle,
$$X_I=G\times S_I \to V_I=G\times_{G_I} S_I.$$
$X_I$ is the pull-back of the $G_I$-bundle $G\to G/G_I$, and 
in particular carries a $G$-invariant connection $\theta_I$ obtained by 
pull-back of the unique  $G$-invariant connection on that bundle. 
For $I\supset J$ there are natural $G$-equivariant 
inclusions $f_I^J:\,X_I\to X_J$, and these are compatible as 
in Section \ref{sec:glue}. The homomorphisms $\varrho_j:\,\pi_1(G_j)\to \U(1)$ 
define flat,  $G$-equivariant bundle gerbes 
$\ca{G}_j=(X_j,L_j,t_j)$ over $V_j$. 

The quotient of the two gerbes on $V_{ij}$, obtained by pulling back 
$\ca{G}_i,\ca{G}_j$ to $X_{ij}$, is just the gerbe defined by 
the homomorphism $\varrho_{ij}:\,\pi_1(G_{ij})\to \U(1)$.
By Lemma \ref{lem:diff} and Proposition \ref{prop:central}(b), 
it follows that this quotient gerbe has a distinguished, 
 equivariant pseudo-line bundle $(E_{ij},s_{ij})$
(where $E_{ij}$ is trivial), with connection $\nabla^{E_{ij}}$ 
induced from the connection $\theta_{ij}$. 
{}From the definition of $\theta_{ij}$, it 
follows that the equivariant error 2-form for this connection is 
the pull-back of the equivariant symplectic form on the 
coadjoint orbit through $\mu_j-\mu_i$.  

We now modify the bundle gerbe connection by adding the equivariant
2-form $(\varpi_j)_G\in\Om^2_G(V_j)$ to the gerbe connection.
Proposition \ref{prop:geom}(d) shows that the equivariant error 2-form
of $\nabla^{E_{ij}}$ with respect to the new gerbe connection
vanishes. The other conditions from the gluing construction in
\ref{sec:glue} are trivially satisfied.  Since the equivariant
3-curvature for the new gerbe connection on $\G_j$ is
$\d_G(\varpi_j)_G=\eta_G|_{V_j}$, we have constructed an equivariant
bundle gerbe with connection, with equivariant curvature-form
$\eta_G$.

\begin{remark}
For $G=\SU(d+1)$ this construction reduces to the the construction in 
terms of transition line bundles: All $L_i,t_i,E_{ij},u_{ijk}$ 
are trivial in this case, hence the entire information on the gerbe 
resides in the functions $s_{ij}:\,(X_{ij})^{[2]}\to \U(1)$ defined 
by the differences $\mu_j-\mu_i$. The condition $\delta s_{ij}=1$ 
for these functions means that $s_{ij}$ defines a line bundle 
$L_{ij}$ over $V_{ij}$, as remarked at the beginning of Section 
\ref{subsec:murray}. The condition $s_{ij}s_{jk}s_{ki}=1$ 
over $X_{ijk}$ is the compatibility condition over triple intersections.
\end{remark}

\section{Pre-quantization of conjugacy classes}
\label{sec:conj}
It is a well-known fact from symplectic geometry  
that a coadjoint orbit $\ca{O}=G.\mu$ through
$\mu\in\t^*_+$ has integral symplectic form, i.e. admits a pre-quantum
line bundle, if and only if $\mu$ is in the weight lattice
$\Lambda^*$. The analogous question for conjugacy classes reads: For
which $\mu\in \Alc$ and $m\in \N$ does the pull-back of the $m$th
power of the basic gerbe $\ca{G}^m$ to the conjugacy class 
$\Co=G.\exp(\mu)$ admit a
pseudo-line bundle, with $m\om_\Co$ as its error 2-form?  For any
positive integer $m>0$ let 
$$\Lambda^*_m=\Lambda^*\cap m\Alc$$ 
be the set of level $m$ weights. As is well-known \cite{pr:lo}, the set
$\Lambda^*_m$ parametrizes the positive energy representations of the
loop group $LG$ at level $m$. 

\begin{theorem}
The restriction of $\G^m$ to a conjugacy class $\Co$ admits a
pseudo-line bundle $\ca{L}$ with connection, with error 2-form
$m\om_\Co$, if and only if $\Co=G.\exp(\mu/m)$ with
$\mu\in\Lambda^*_m$. Moreover $\ca{L}$ has an equivariant extension in
this case, with $m\om_\Co$ as its {\em equivariant} error 2-form .
\end{theorem} 

\begin{proof}
Given a conjugacy class $\Co\subset G$, let 
$\mu\in m\Alc$ be the unique point with $g:=\exp(\mu/m)\in\Co$,
and let $K=G_g$ so that $\Co=G/K$. Pick an index $j$ with $\Co\subset
V_j$, and let
$$\nu=m\Psi_j(g)=\mu-m\mu_j.$$
Then 
$$ G_\mu\subset K\subset G_\nu.$$
Let $\O_\mu,\O_\nu\subset\g$ denote the adjoint orbits of $\mu,\nu$,
and $(\om_\mu)_G,(\om_\nu)_G$ their equivariant symplectic forms.  The
pull-back $\iota_\Co^*\ca{G}^m$ is the gerbe over $G/K$ defined 
as in Section \ref{sec:prin}
by
the homomorphism $\varrho\in \Hom(\pi_1(K),\U(1))$, given as a
composition
$$\pi_1(K)\to \pi_1(G_j)\to \U(1),$$
where the fist map is push-forward under the inclusion $K\hra G_j$,
and the second map is the homomorphism defined by the element
$m\mu_j\in\t$ for $G_j$.

Suppose now that $\mu\in\Lambda^*_m$. Then $m\mu_j$ equals
$-\nu$ up to a weight lattice vector, which means that $\varrho$ 
is the image of $-\nu\in (\k^*)^K$ 
in the exact sequence \eqref{eq:sequence}. 
Hence, Proposition \ref{prop:central} says that we 
we obtain an 
equivariant pseudo-line bundle for $\iota_\Co^*\ca{G}^m$, with 
equivariant error 2-form 
$$ \Psi_j^*(\om_{\nu})_G-m\,\iota_\Co^*(\varpi_j)_G=m\,\om_\Co.$$
Here we have used part (b) of Proposition \ref{prop:geom}.

Conversely, suppose that $\G^m|_\Co$ admits a pseudo-line bundle with
error 2-form $m\,\om_\Co$. Consider
the pull-back of $\G$ under the exponential map $\exp:\,\g\to G$. 
The pull-back
$\exp^*\eta\in\Om^3(\g)$ is exact, and the homotopy operator for the
linear retraction of $\g$ to the origin defines a 2-form
$\varpi\in\Om^2(\g)$ with $\d\varpi=\exp^*\eta$. As in Proposition
\ref{prop:geom}, one shows that for any adjoint orbit $\O\subset \g$,
with $\exp\O=\Co$,
$$ \iota_\O^*\varpi=\exp^*\om_{\Co}-\om_\O$$
where $\om_\O$ is the symplectic form on $\O$. 
In particular this applies to $\O=\O_{\mu/m}$.
Choose a pseudo-line bundle for $\exp^*\G$ with error 2-form 
$-\varpi$. We then have two pseudo-line bundles for $\exp^*\G^m|_{\O}$ 
obtained by restricting the $m$th power of the pseudo-line bundle
for $\exp^*\G$ or by pulling back the pseudo-line bundle for 
$\Co$. Their quotient is a line bundle over $\O$, with 
curvature the difference of the error 2-forms:
$$ m(\exp^*\om_{\Co}-\iota_{\O_\mu}^*\varpi)=m\om_\O.$$
Thus $m(\mu/m)=\mu$ must be in the weight lattice. 
\end{proof}

\begin{remark}
Z. Shahbazi has proved that if $\G$ is a gerbe with connection
over a manifold $M$, with curvature 3-form $\eta$, and $\Phi:\,N\to M$ 
is a map with $\Phi^*\eta+\d\om=0$, then the pull-back gerbe $\Phi^*\G$ admits 
a pseudo-line bundle, with $\om$ as its error 2-form, if and only if 
the pair $(\eta,\om)$ defines an integral element of the relative 
de Rham cohomology $H^3(\Phi,\R)$. This means that for any smooth 2-cycle
$S\subset N$, and any smooth 3-chain $B\subset M$ with boundary $\Phi(S)$, 
one must have $\int_B\eta-\int_S\om\in\Z$. The particular case
where the target of $\Phi$ is a Lie group $G$ is relevant for the 
pre-quantization of group-valued moment maps \cite{al:mom}.
\end{remark}

\begin{appendix}

\section{Proof of Lemma \ref{lem:covers}}\label{sec:cover}
In this Appendix we prove Lemma \ref{lem:covers}, concerning the 
construction of a certain cover $U_I$ of $M$ from a given cover 
$V_j$. Write $M=\coprod_I A_I$ where 
$$ A_I=\bigcap_{i\in I}V_i\backslash \bigcup_{j\not\in I}V_j.$$
Notice that $\ol{A_I}\subset \bigcup_{J\subset I} A_J.$ 
By induction on the cardinality $k=|I|$ we will 
construct open sets $U_J\subset V_J$, 
having the following properties:  
\begin{enumerate}
\item
the closure 
$\ol{U_I}$ does not meet $\ol{U_J}$ for $|J|\le |I|$ unless $J\subset I$, 
\item 
each $\ol{A_I}$ is contained in the union of $U_J$ with 
$J\subset I$.
\end{enumerate}
The induction starts at $k=0$, taking $U_\emptyset=\emptyset$. Suppose 
we have constructed open sets $U_I$ with $\ol{U_I}\subset V_I$ for $|I|<k$, such that 
the properties (a),(b) hold for all $|I|<k$. For $|I|=k$ consider the 
subsets 
$$ B_I:=A_I\backslash \Big(\bigcup_{{ J\subset I},|J|<k}U_J\Big).$$
Note that (unlike $A_I$) the set $B_I$ is closed.
$B_I$ does not meet $\ol{A_J}$ unless $I\subset J$, and 
it also does not meet  $\ol{U_J}$ for $|J|<k$ unless $J\subset I$. 
That is, $B_I$ is disjoint from 
$$ C_I:=
\bigcup_{{J\not\subset I},|J|<k}\ol{U_J}
\cup \bigcup_{{K\not\supset I}}\ol{A_K}
$$
Choose open sets $U_I$ for $|I|=k$ with $B_I\subset U_I\subset
\ol{U_I}\subset M\backslash C_I$, and such that the closures of the
sets $U_I$ for distinct $I$ with $|I|=k$ are disjoint.  The new
collection of subsets will satisfy the properties (a),(b) for $|I|\le
k$. We next show that $V_i'=M\backslash \bigcup_{J\not\ni i}\ol{U_J}$
is a cover of $M$. Write $M=\coprod_I D_I$ with
$D_I=\ol{U_I}\backslash \bigcup_{|J|<|I|}\ol{U_J}$.  Then $D_I\cap
\ol{U}_J=\emptyset$ unless $I\subset J$, so $D_I$ is contained in each
$V_i'$ with $i\in I$.  In particular $\bigcup_i V_i'=M$. Finally
$\ol{V_i'}\subset \bigcup_{I\ni i}\ol{U_I}\subset V_i$.  This
completes the proof of Lemma \ref{lem:covers}. Note that if the $V_i$
were invariant under an action of a compact group $G$, the $U_I$ could
be taken $G$-invariant also.
\end{appendix}

\bibliographystyle{amsplain}

\def\polhk#1{\setbox0=\hbox{#1}{\ooalign{\hidewidth
  \lower1.5ex\hbox{`}\hidewidth\crcr\unhbox0}}} \def\cprime{$'$}
  \def\cprime{$'$} \def\polhk#1{\setbox0=\hbox{#1}{\ooalign{\hidewidth
  \lower1.5ex\hbox{`}\hidewidth\crcr\unhbox0}}} \def\cprime{$'$}
\providecommand{\bysame}{\leavevmode\hbox to3em{\hrulefill}\thinspace}
\providecommand{\MR}{\relax\ifhmode\unskip\space\fi MR }
\providecommand{\MRhref}[2]{%
  \href{http://www.ams.org/mathscinet-getitem?mr=#1}{#2}
}
\providecommand{\href}[2]{#2}

\end{document}